\documentclass{article}
\def\Z{\mathbb Z}
\def\N{\mathbb N}
\def\A{\mathcal A}

\def\pf{\begin{proof}}
\def\pfk{\end{proof}}
\def\AC{\mathrm{AC}}

\usepackage{amssymb,amsmath,amsthm}

\newtheorem{lem}{Lemma}[section]
\newtheorem{obs}[lem]{Observation}
\newtheorem{prop}[lem]{Proposition}
\newtheorem{thm}[lem]{Theorem}
\newtheorem{coro}[lem]{Corollary}

\theoremstyle{definition}
\newtheorem{de}[lem]{Definition}
\newtheorem{pozn}[lem]{Remark}

\begin{document}
\title{Balances and Abelian Complexity of a Certain Class of Infinite Ternary Words}
\author{Ond\v rej Turek}


\maketitle

\begin{center}
Laboratory of Physics, Kochi University of Technology\\
Tosa Yamada, Kochi 782-8502, Japan\\
email: \texttt{ondrej.turek@kochi-tech.ac.jp}
\end{center}
%
%
\bigskip
\begin{abstract}
A word $u$ defined over an alphabet $\mathcal{A}$ is $c$-balanced ($c\in\mathbb{N}$) if for all pairs of factors $v$, $w$ of $u$ of the same length and for all letters $a\in\mathcal{A}$, the difference between the number of letters $a$ in $v$ and $w$ is less or equal to $c$. In this paper we consider a ternary alphabet $\mathcal{A}=\{L,S,M\}$ and a class of substitutions $\varphi_p$ defined by $\varphi_p(L)=L^pS$, $\varphi_p(S)=M$, $\varphi_p(M)=L^{p-1}S$ where $p>1$.
We prove that the fixed point of $\varphi_p$, formally written as $\varphi_p^\infty(L)$, is 3-balanced and that its Abelian complexity is bounded above by the value 7, regardless of the value of $p$. We also show that both these bounds are optimal, i.e. they cannot be improved.
\end{abstract}
%
%
%
\maketitle
\section*{Introduction}

The balance property is a notion connected with Sturmian words
from very first beginning of their investigation. In
\cite{morse2}, Sturmian words were defined as aperiodic words
with the smallest possible factor complexity. Already in the same
article, Hedlund and Morse observed that Sturmian words show also
the smallest discrepancy in occurrences of letters. To quote
precisely their result, let us denote by $|w|$ the length of the
word $w$  and by $|w|_a$ the number of occurrences of letter $a$
in $w$. Hedlund and Morse proved that an infinite aperiodic word
$u$ over the alphabet $\{0,1\}$ is Sturmian if and only if for all
pairs $w$, $v$ of factors of $u$ with $|w|=|v|$ it holds
 $||w|_0-|v|_0| \leq 1$. Let us note that the  letter  $0$ is not
preferred as   in binary alphabet the relations $|w|=|v|$ and
$||w|_0-|v|_0| \leq 1$ imply  the inequality $||w|_1-|v|_1| \leq
1$ as well.

During past 70 years many other characterizations of Sturmian
words have appeared, for their overview see \cite{lothaire}. Each
of these characterizations may serve and serves for generalization
of Sturmian words to multiliteral alphabets, cf. \cite{BaPeStar}. Nevertheless, the
balance property seems to be the most complicated to deal with.
Only a few results are known about words satisfying the so-called
$c$-balanced property.

 Let us recall that  an infinite word $u$
over an alphabet $\mathcal{A}$ is $c$-balanced if for all letters
$a\in \mathcal{A}$ and all pairs  $v, w$ of factors of $u$ with
$|v|=|w|$ it holds $||v|_a-|w|_a| \leq c$. 
Note that Sturmian words are $1$-balanced in this terminology. The set of  $c$-balanced words
differs substantially from all other generalizations of Sturmian
words. Neither generic  Arnoux-Rauzy  word nor generic  word
coding interval exchange transformation are $c$-balanced, see
\cite{CaFeZa} and  \cite{Ad2002}.

 In \cite{adamczewski},
Adamczewski studies whether fixed point of a primitive
substitution  is $c$-balanced for some constant $c$. He shows that
the existence of such $c$ depends only on the spectrum  of the incidence
matrix of the substitution.  However, the minimal value of $c$ cannot
be deduced from the spectrum.  In \cite{turek} and \cite{BPT}, the
minimal value of $c$ is determined for binary fixed points of
canonical substitutions associated with quadratic Pisot numbers.
 The notion ``canonical substitution associated
with a number $\beta >1$'' comes from  positional  numeration
systems with the base $\beta$, see \cite{fabre}. Generally speaking, it is a very complicated problem to determine minimal
value $c$ for a ternary balanced word, let alone for words over alphabets of higher cardinalities. Despite
a common belief that the Tribonacci word is $2$-balanced, the first
proof of this fact has appeared just one year ago in
\cite{RSZ} (the Tribonacci word is the fixed point of the
substitution  $A\mapsto AB$, $B\mapsto AC$, $C\mapsto A$). In this article we
provide minimal value of $c$ for a certain class of ternary words,
namely for fixed points of substitutions

\begin{equation}\label{nase}
L\mapsto L^pS\,, \quad S\mapsto M\,, \quad M\mapsto L^{p-1}S
\end{equation}
 with
the parameter $p>1$. These substitutions are canonical substitutions
associated with cubic Pisot numbers $\beta>1$, roots of
polynomials $x^3 -px^2-x+1$ (cf. \cite{KP}). Let us recall that the
Tribonacci substitution is associated with a numeration system as well.

The definition of a $1$-balanced  word may be reformulated
equivalently using Parikh vectors.
Inspired by this fact, Richome, Saari and Zamboni introduced
the Abelian complexity $\AC(n)$ of infinite word. In their notation,
Sturmian words are aperiodic words with $\AC(n) = 2$ for all $n\in
\mathbb{N}$. The
 question on existence of
 words with constant Abelian complexity  is natural. It was shown in \cite{CR}
 that for  $k\geq 4$ no words with $\AC(n) = k$ exist. On the other hand,
 words with  $\AC(n) = 3$ can be found in \cite{RSZ2}.

The relation between Abelian complexity and balance property is
not straightforward.  It is easy
 to see that an infinite word $u$ is balanced  if and only if its Abelian complexity is bounded.
 Moreover, if the Abelian complexity of $u$ is bounded by $k$, then $u$ is $k-1$
 balanced. The Tribonacci case  shows that the opposite implication
 is not valid: according to \cite{RSZ}, the Abelian complexity of the Tribonacci word
takes all values in the set $\{3,4,5,6,7\}$. The fixed point of the substitutions
\eqref{nase}  studied in this article has the same property.

\section{Preliminaries}

Let $\A$ be a finite alphabet. A concatenation of letters in $\A$
is called a \emph{word}. The set $\A^*$ of all finite words over $\A$ equipped with the empty
word $\epsilon$ and the operation of concatenation is a free monoid.
The \emph{length} of the word $w\in\A^*$, denoted by $|w|$, represents the number of its letters.

One may also consider infinite words $u=u_0u_1u_2\cdots$; the set of infinite words
over the alphabet $\A$ is denoted by $\A^\N$.

A word $w$ is called a \emph{factor} of $v\in\A^*$ or $\A^\N$
if there exist words $w^{(1)}\in\A^*$ and $w^{(2)}\in\A^*$ or
$w^{(2)}\in\A^\N$, respectively, such that $v=w^{(1)}ww^{(2)}$. The word $w$
is called a \emph{prefix} of $v$, if $w^{(1)}=\epsilon$. It is a \emph{suffix} of
$v$, if $w^{(2)}=\epsilon$.

Let $w\in\A^{\N}$. For $k\in\N$, the symbol $w^k$ denotes the concatenation $\underbrace{ww\cdots w}_{k\,\text{times}}$. We set $w^0=\epsilon$. Let a word $v\in\A^*$ have the prefix $w^k$, $k\in\N$. Then the symbol $w^{-k}v$ denotes the word satisfying $w^kw^{-k}v=v$. Similarly, if a word $v\in\A^\N$ has the suffix $w^k$ for a $k\in\N$, then $vw^{-k}$ denotes the word with the property $vw^{-k}w^k=v$.

A \emph{morphism} on the free monoid $\A^*$ is a map $\varphi:\A^*\to\A^*$
satisfying $\varphi(vw)=\varphi(v)\varphi(w)$
for all $v,w\in\A^*$. Obviously, the morphism $\varphi$
is determined if we define $\varphi(a)$ for all $a\in\A$.

A morphism $\varphi$ is called a \emph{substitution}, if $\varphi(a)
\neq\epsilon$ for all $a\in\A$ and if there is an $a'\in\A$
such that $|\varphi(a')|>1$. An infinite word $u$ is said to be
a \emph{fixed point} of the substitution $\varphi$, or invariant under
the substitution $\varphi$, if
\begin{equation}\label{eq:word}
\varphi(u_0)\varphi(u_1)\varphi(u_2)\cdots = u_0u_1u_2\cdots\,.
\end{equation}
If we naturally extend the action of $\varphi$ to infinite words, we may rewrite \eqref{eq:word} simply as $\varphi(u)=u$.

\subsection{Balance properties}


An infinite word $u$ is $c$-balanced, if for every
$a\in\A$ and for every pair of factors $v$, $w$ of $u$ such that $|v|=|w|$, it holds
$\left||v|_a-|w|_a\right|\leq c$. This property determines the discrepancy of occurrences of letters in the word $u$. However, it turns out that if the cardinality of $\A$ is higher than two, it is useful to have more detailed information, namely what is the discrepancy of occurrences of each particular letter. For this purpose we introduce the following notion:

\begin{de}
Let $u$ be an infinite word over the alphabet $\A$ and let $a\in\A$. The word $u$ is said to be \emph{$c$-balanced with respect to the letter $a$}, if
$$
\left|\, |v|_a-|w|_a \,\right|\leq c
$$
for all pairs of factors $v$, $w$ of $u$ of the same length.
\end{de}

\subsection{Abelian complexity}\label{PrelAC}

Let us consider an alphabet $\A$ with $k$ elements, i.e. $\A=\{a_1,\ldots,a_k\}$, and an infinite word $u$ over $\A$. For any factor $w$ of $u$, its \emph{Parikh vector} is the $k$-tuple $\Psi(w)=(|w|_{a_1},\ldots,|w|_{a_k})$. Let the symbol $\mathcal{F}_u(n)$ denote the set of all factors of $u$ of the length $n$. Then the \emph{Abelian complexity} of the word $u$ is a function $\mathrm{AC}:\N\to\N$ defined by
\begin{equation}\label{AC}
\mathrm{AC}(n)=\#\left\{\Psi(w)\,\left|\,w\in\mathcal{F}_u(n)\right.\,\right\}\,.
\end{equation}
On the right hand side of \eqref{AC} there is the cardinality of the set of Parikh vectors of all factors of $u$ of the length $n$. In the sequel we will denote this set by $\mathcal{P}_u(n)$, i.e.
$$
\mathcal{P}_u(n)=\left\{\Psi(w)\,\left|\,w\in\mathcal{F}_u(n)\right.\,\right\}\,.
$$

\subsection{On the word studied in this paper}

From now on, we will focus on a special class of substitutions on the ternary alphabet $\{L,S,M\}$. For any integer $p>1$, we denote by $\varphi_p$ the substitution
given by
 \begin{equation}\label{varphi}
 \begin{array}{rcl}
 \varphi_p(L) &=& L^pS\\
 \varphi_p(S) &=& M\\
 \varphi_p(M) &=& L^{p-1}S
 \end{array}
 \end{equation}
The substitution $\varphi_p$ has a unique fixed point, namely
$$
u^{(p)}=\lim_{n\to\infty}\varphi_p^n(L).
$$

If the results of \cite{adamczewski} are applied on $u^{(p)}$, one finds out that there is a constant $c$ such that the word $u^{(p)}$ is $c$-balanced, but it is not known what the value of $c$ is and how it depends on $p$. This is the main aim of this paper -- to determine $c$.

The fact that $u^{(p)}$ is balanced immediately implies that the Abelian complexity function of $u^{(p)}$ is bounded, see Introduction. The second aim of this paper is thus to find the optimal bound for $\mathrm{AC}(n)$.

\begin{pozn}
The elements of $\A$ are usually denoted by numbers: $0,1,2$ etc. We have considered this notation, but we believe that the paper becomes more transparent if letters are used. The choice of $L$, $S$ and $M$ has its roots in the fact that the word $u^{(p)}$ is a fixed point of a substitution associated with a number $\beta>1$, cf. Introduction. Let $\Z_\beta$ denote the set of numbers which can be written in the form $x=x_k\beta^k + \cdots +x_1\beta + x_0$ for non-negative integers $x_j$. It can be shown (cf. \cite{thurston}) that when the elements of $\Z_\beta$ are drawn on the real line, there are exactly three types of distances between neighbouring points. If we assign the letters $L$, $M$ and $S$ to the longest, the medium and the shortest distance, respectively, then the order of distances on the real line corresponds exactly to the order of the letters $L,S,M$ in the infinite word $u^{(p)}$.
\end{pozn}

\section{Main result and the proof outline}

We begin by the formulation of the main result of the paper.

\begin{thm}\label{balance}
Let $u^{(p)}$ be the infinite word invariant under the morphism $\varphi_p$ given by \eqref{varphi}. Then $u^{(p)}$ is
\begin{itemize}
\item $3$-balanced with respect to the letter $L$,
\item $2$-balanced with respect to the letter $S$,
\item $2$-balanced with respect to the letter $M$,
\end{itemize}
and none of these bounds can be improved.
\end{thm}
The theorem has the following trivial consequence:

\begin{coro}
The infinite word $u^{(p)}$ is $3$-balanced and this bound is optimal, i.e. it cannot be improved.
\end{coro}

Since the proof is long and slightly complicated, we will split it into four sections and proceed in the following way:
\begin{itemize}
\item[1.] We prove that $u^{(p)}$ is $2$-balanced with respect to the letter $M$.
\item[2.] We prove that $u^{(p)}$ is $2$-balanced with respect to the letter $S$.
\item[3.] We prove that $u^{(p)}$ is $3$-balanced with respect to the letter $L$.
\item[4.] We show that none of the bounds can be improved.
\end{itemize}

\section{Properties of the word $u^{(p)}$}

As we have explained in Preliminaries, the word $u^{(p)}$ is a fixed point of $\varphi_p$, i.e.
$$
u^{(p)}=\varphi_p(u)=\varphi_p(u_0)\varphi_p(u_1)\varphi_p(u_2)\cdots
$$
In this sense each letter of $u^{(p)}$ can be regarded as the image, or a factor of the image, of another letter of $u^{(p)}$. In view of the definition of $\varphi_p$, cf. \eqref{varphi}, each segment $\varphi_p(u_j)$ has the structure $L^kY$ for $k\in\{0,p-1,p\}$ and $Y\neq L$. The letters $S$ and $M$ are thus ``terminating symbols'' which cut $u^{(p)}$ to images of individual letters.
This fact is particularly important when a factor $v$ of $u^{(p)}$ is given and one needs to find a factor $x$ of $u^{(p)}$ such that $\varphi_p(x)=v$. It holds:

\begin{obs}\label{JednozVzor}
Let $vY$ be a factor of $u^{(p)}$ such that $Y\in\{S,M\}$ and let one of the following conditions be satisfied:
\begin{itemize}
\item[\textit{(i)}] The first letter of $v$ is $M$,
\item[\textit{(ii)}]  $v$ has the prefix $L^p$,
\item[\textit{(iii)}]  $SvY$ or $MvY$ is a factor of $u^{(p)}$.
\end{itemize}
Then there is a unique factor $x$ of $u^{(p)}$ satisfying $\varphi_p(x)=vY$.\\
\end{obs}

\pf
Any of the conditions (i), (ii) and (iii) together with $Y\in\{S,M\}$ ensures that $vY$ is an image of certain factor $x$, and it is obvious from the definition of $\varphi_p$ that $\varphi_p(x)=\varphi_p(y)\Rightarrow x=y$.
\pfk

\begin{obs}\label{pripustneVzor}
\textit{(i)}\quad Let $X$ be a letter occuring in $\varphi_p(u_j)$ for a $j\in\N_0$. Then:
\begin{itemize}
\item If $X=M$, then $u_j=S$ and $\varphi_p(u_j)=X$.
\item If $X=S$, then either $u_j=L$ and $\varphi_p(u_j)=L^pX$, or $u_j=M$ and $\varphi_p(u_j)=L^{p-1}X$.
\end{itemize}
\textit{(ii)}\quad If $XL^kY$ is a factor of $u^{(p)}$ and $X\neq L$, $Y\neq L$, $k\neq0$, then $Y=S$ and $k\in\{p-1,p\}$.
\end{obs}

The following observation describes the possible neighbours of each of the letters $L,S,M$ in the word $u^{(p)}$.

\begin{obs}\label{pripustne}
The sequence of letters in the word $u^{(p)}$ conform to these rules:
\begin{itemize}
\item[\textit{(i)}] Each letter $S$ in $u^{(p)}$ is preceded by $L$ and followed either by $L$ or by $M$.
\item[\textit{(ii)}]  Each letter $M$ in $u^{(p)}$ is preceded by $S$ and followed by $L$.
\end{itemize}
\end{obs}

\pf
\textit{(i)} Each $S$ is the last letter of $\varphi_p(u_j)$ for $u_j=L$ or $u_j=M$ according to Observation \ref{pripustneVzor}, i.e. it is the last letter of the block $L^pS$ or $L^{p-1}S$, thus is preceded by $L$. The letter $S$ is followed by the first letter of $\varphi_p(u_{j+1})$, which can be either $L$ or $M$ (cf. the substitution rule $\varphi_p$).\\
\textit{(ii)} Each $M$ is equal to $\varphi_p(u_j)$ for $u_j=S$. We already know from (i) that $u_{j-1}=L$ and $u_{j+1}\in\{L,M\}$, therefore the $M$ is preceded by the last letter of $\varphi_p(L)$ (which is $S$) and followed by the first letter of $\varphi_p(L)$ or $\varphi_p(M)$ (which is $L$).
\pfk

In order to understand the structure of the word $u^{(p)}$, it is useful to describe possible segments $z$ in factors of $u^{(p)}$ of the type $SzS$ and $MzM$. This is done in the next two observations.

\begin{obs}\label{SzS}
Let $Sz'S$ be a factor of $u^{(p)}$ such that $|z'|_S=0$. Then one of the following equalities holds:
\begin{itemize}
\item $z'=L^p$,
\item $z'=ML^p$,
\item $z'=ML^{p-1}$.
\end{itemize}
\end{obs}

\pf
Observations \ref{JednozVzor} and \ref{pripustneVzor} imply that $z'S=\varphi_p(\check{z}')$, where $\check{z}'$ is a factor of $u^{(p)}$ the last letter of which is either $L$ or $M$ and which is preceded in $u^{(p)}$ by either $L$ or $M$.

Since $\varphi_p(\check{z}')$ contains only one $S$, namely its last letter, all letters of $\check{z}'$ except the last one have to be different from $L$ and $M$. Therefore $\check{z}'=S^kX$, where $X\in\{L,M\}$ and $k\geq0$.

Taking into account Observation \ref{pripustne}, we infer that $k=0$ or $k=1$. Therefore only four situations are possible: $\check{z}'=L$, $\check{z}'=M$, $\check{z}'=SL$, $\check{z}'=SM$. Moreover, since $\check{z}'$ is preceded by either $L$ or $M$, it cannot hold $\check{z}'=M$, cf. Observation~\ref{pripustne}. Therefore $\check{z}'$ is equal to one of the factors $L$, $SL$, $SM$, which implies that $z$ is equal to one of the factors $L^p$, $ML^p$, $ML^{p-1}$.
\pfk

\begin{obs}\label{MzM}
Let $MzM$ be a factor of $u^{(p)}$ such that $|z|_M=0$. Then one of the following equalities holds:
\begin{itemize}
\item $z=(L^pS)^p$,
\item $z=L^{p-1}S(L^pS)^p$,
\item $z=L^{p-1}S(L^pS)^{p-1}$.
\end{itemize}
Consequently, $p^2+p-1\leq |z|\leq p^2+2p$.
\end{obs}

\pf
It follows from Observations \ref{JednozVzor} and \ref{pripustneVzor} that $MzM=\varphi_p(Sz'S)$, where $Sz'S$ is a factor of $u^{(p)}$ such that $|z'|_S=0$. Therefore $z=\varphi_p(z')$ and, according to Observation \ref{SzS}, $z'\in\{\varphi_p(L^p), \varphi_p(ML^p), \varphi_p(ML^{p-1})\}$.
\pfk

Many times we will need to compare the number of letters $L,S,M$ in a factor of $u^{(p)}$ and in its image. The substitution rule \eqref{varphi} leads to the equalities
$$
|\varphi_p(v)|_L=p|v|_L+(p-1)|v|_M\,,\quad |\varphi_p(v)|_S=|v|_L+|v|_M\,,\quad |\varphi_p(v)|_M=|v|_S\,,
$$
which can be inverted subsequently:

\begin{prop}\label{vzor}
For any factor $v$ of $u^{(p)}$ it holds
\begin{align*}
|v|_L&=|\varphi_p(v)|_L-(p-1)|\varphi_p(v)|_S\,, \\
|v|_S&=|\varphi_p(v)|_M\,, \\
|v|_M&=-|\varphi_p(v)|_L+p|\varphi_p(v)|_S\,, \\
|v|&=|\varphi_p(v)|_S+|\varphi_p(v)|_M\,.
\end{align*}
\end{prop}

\section{Balance bound with respect to the letter $M$}
We begin the proof of Theorem \ref{balance} by its second statement, i.e. we show at first that $u^{(p)}$ is 2-balanced with respect to the letter $M$. As we will see, the determination of the balance bound with respect to the letter $M$ is by far the most complicated part of the work.

\begin{thm}\label{balanceM}
Let $v$, $w$ be factors of $u^{(p)}$ such that $|v|=|w|$. Then
$$
\left|\, |v|_M-|w|_M\, \right|\leq2\,.
$$
\end{thm}

\subsection*{Proof of Theorem \ref{balanceM}}

We will proceed by contradiction. Let us assume that there exist factors $v$, $w$ of $u^{(p)}$ such that $|v|=|w|=n$ and
\begin{equation}\label{unbalance M}
|v|_M-|w|_M>2\,.
\end{equation}
Let $n$ be the \emph{minimal} number with this property.

We denote $v=v_1\cdots v_n$, $w=w_1\cdots w_n$. The minimality of $n$ implies
\begin{gather}
v_1=M\,, \quad v_n=M\,, \label{vMM} \\
w_1\neq M\,, \quad w_n\neq M\,, \label{wxMxM}\\
|v|_M-|w|_M=3\,. \label{vw3M}
\end{gather}

\subsubsection*{Stage 1: Introduction of $f$, $g$}

If we apply Observations \ref{JednozVzor} and \ref{pripustneVzor}, Eq. \eqref{vMM} implies that the factor $v$ is an image of certain factor of $u^{(p)}$ whose first and last letters equal $S$. This allows us to define a factor $f$ of $u^{(p)}$ in this way:
\begin{equation}\label{f}
\varphi_p(SfS)=v\,.
\end{equation}

The factor $w$ is not ready for a direct application of Observation \ref{JednozVzor} because of \eqref{wxMxM}. For that reason we at first extend the factor $w$ to both sides up to the closest letter $M$, i.e. we put
\begin{equation}\label{defw'}
w'=Mw^{(1)}ww^{(2)}M\,,
\end{equation}
where $w'$ is a factor of $u^{(p)}$ and $|w^{(1)}|_M=|w^{(2)}|_M=0$. Now we can define a factor $g$ of $u^{(p)}$ by the relation
\begin{equation}\label{g}
\varphi_p(SgS)=w'\,.
\end{equation}
Let us show that the factor $g$ is shorter than $v$ and $w$:

\begin{prop}\label{g kratsi}
It holds $|g|\leq n-2(p^2-1)$.
\end{prop}
\pf
Since $|v|_M=|w|_M+3$, the factor $v$ contains at least 3 letters $M$, and thus $v=M\cdots M\cdots M$. Observation \ref{MzM} then implies that $|v|\geq 1+(p^2+p-1)+1+(p^2+p-1)+1=2p^2+2p+1$.

Since $|w|=|v|$, it holds $|w|\geq 2p^2+2p+1$, hence necessarily $|w|_M\geq1$ according to Observation \ref{MzM}. But then $|v|_M\geq 4$ and
\begin{equation}\label{vMMM}
v=Mz^{(1)}M\cdots Mz^{(2)}M\,,
\end{equation}
where $|z^{(1)}|_M=|z^{(2)}|_M=0$ and $|z^{(j)}|\geq p^2+p-1$ for $j=1,2$.

One more application of Observation \ref{MzM} gives $|v|\geq 1+(p^2+p-1)+1+(p^2+p-1)+1+(p^2+p-1)+1=3p^2+3p+1$, hence $|w|\geq 3p^2+3p+1$ and necessarily $|w|_M\geq2$, i.e.
\begin{equation}\label{wMMM}
w=\hat{w}^{(1)}M\hat{w} M\hat{w}^{(2)}\,,
\end{equation}
where $|\hat{w}^{(1)}|_M=|\hat{w}^{(2)}|_M=0$.

We deduce from Eqs. \eqref{vMMM}, \eqref{wMMM}, from $|z^{(j)}|\geq p^2+p-1$ and from the minimality of $n$ that $|\hat{w}^{(j)}|\geq p^2+p$, $j=1,2$.

The factor $w'$ (cf. \eqref{defw'} and \eqref{wMMM}) is given by $w'=Mw^{(1)}\hat{w}^{(1)}M\hat{w} M\hat{w}^{(2)}w^{(2)}M$, hence Equality \eqref{g} together with Observation \ref{pripustneVzor} imply
$$
SgS=Sg^{(1)}S\hat{g} Sg^{(2)}S\,,
$$
where $\varphi_p(g^{(1)})=w^{(1)}\hat{w}^{(1)}$, $\varphi_p(g^{(2)})=\hat{w}^{(2)}w^{(2)}$ and $\varphi_p(S\hat{g}S)=M\hat{w}M$. Hence $|g^{(1)}|_S=|g^{(2)}|_S=0$, and consequently, using Observation \ref{SzS}, $|g^{(1)}|\leq p+1$ and $|g^{(2)}|\leq p+1$.
The sought inequality $|g|\leq n-2(p^2-1)$ follows from these relations:
\begin{gather*}
n=|w|=|\hat{w}^{(1)}\varphi_p(S\hat{g}S)\hat{w}^{(2)}|=|\hat{w}^{(1)}|+|\varphi_p(S\hat{g}S)|+|\hat{w}^{(2)}|\geq 2(p^2+p)+|S\hat{g}S|\,,
\end{gather*}
$\;|g|=|g^{(1)}S\hat{g}Sg^{(2)}|=|g^{(1)}|+|S\hat{g}S|+|g^{(2)}|\leq 2(p+1)+|S\hat{g}S|$.
\pfk


\begin{prop}
Factors $f$, $g$ satisfy
\begin{gather}
|f|_S-|g|_S=1\,, \label{f,g S}\\
|f|_M-|g|_M\geq(p+1)\cdot(|f|-|g|)+2-p\,. \label{*}
\end{gather}
\end{prop}

\pf
Both relations will be proved using the properties of $v$, $w'$.
Let us begin with Eq. \eqref{f,g S}. Considering Proposition \ref{vzor}, we have
\begin{align*}
|v|_M&=|\varphi_p(SfS)|_M=|SfS|_S=|f|_S+2\,,\\
|w'|_M&=|\varphi_p(SgS)|_M=|SgS|_S=|g|_S+2\,.
\end{align*}
Moreover, $|w'|_M=|Mw^{(1)}ww^{(2)}M|_M=|w|_M+2$. Therefore, taking into account Eq. \eqref{vw3M}, we obtain
$$
|f|_S-|g|_S=|v|_M-|w'|_M=|v|_M-|w|_M-2=3-2=1\,.
$$
Now we proceed to Eq. \eqref{*}. It holds $|\varphi_p(z)|=(p+1)|z|_L+|z|_S+p|z|_M$ for any factor $z$ (cf. the substitution rule \eqref{varphi}). The identity $|z|=|z|_L+|z|_S+|z|_M$ allows one to eliminate $|z|_L$, hence
\begin{equation}\label{varphi(z)}
|\varphi_p(z)|=(p+1)|z|-p|z|_S-|z|_M\,.
\end{equation}
Since $|v|=|w|$ and $|w'|=|Mw^{(1)}ww^{(2)}M|\geq |w|+2$, we have $|v|-|w'|\leq-2$, equivalently $|\varphi_p(SfS)|-|\varphi_p(SgS)|\leq-2$. If we apply \eqref{varphi(z)} to the last inequality, we obtain
$$
(p+1)(|f|-|g|)-p(|f|_S-|g|_S)-(|f|_M-|g|_M)\leq-2\,.
$$
Now we substitute here from \eqref{f,g S} which leads to \eqref{*}.
\pfk

There is another useful statement, namely Proposition \ref{f-g<=0}, but to prove it we need one more observation which gives an estimate of the number of occurences of the letter $M$ in a factor of $u^{(p)}$ of a given length:

\begin{obs}\label{obsahM}
\textit{(i)}\quad For every factor $\hat{v}$ of $u^{(p)}$ it holds $|\hat{v}|_M\leq1+\frac{|\hat{v}|-1}{p^2+p}$. \\
\textit{(ii)}\quad If moreover $\hat{v}$ has a prefix or a suffix of the length greater or equal to $\Delta$ that does not contain $M$, it holds $|\hat{v}|_M\leq\left\lceil\frac{|\hat{v}|-\Delta}{p^2+p}\right\rceil$.
\end{obs}
\pf
(i) Observation \ref{MzM} says that if $MzM$ is a factor of $u^{(p)}$, then $|zM|\geq p^2+p$. Consequently $|\hat{v}|_M\leq\left\lceil\frac{|\hat{v}|}{p^2+p}\right\rceil$. The inequality $\left\lceil\frac{|\hat{v}|}{p^2+p}\right\rceil\leq1+\frac{|\hat{v}|-1}{p^2+p}$ holds trivially.\\
(ii) Let $\hat{v}=\hat{v}'\hat{v}''$ and
$$
|\hat{v}'|=\Delta\quad\wedge\quad|\hat{v}'|_M=0\qquad \text{or} \qquad |\hat{v}''|=\Delta\quad\wedge\quad|\hat{v}''|_M=0\,.
$$
Then, employing the result of (i), one has $|\hat{v}|_M=|\hat{v}'|_M+|\hat{v}''|_M\leq0+\left\lceil\frac{|\hat{v}|-\Delta}{p^2+p}\right\rceil$.
\pfk

\begin{prop}\label{f-g<=0}
It holds $|f|-|g|\leq 0$, i.e. the factor $f$ is not longer than $g$.
\end{prop}

\pf
Let us suppose that the contrary is true, i.e. $|f|-|g|=d>0$. Let $f'$ be the suffix of $f$ of the length $d$ and let $\hat{f}=f(f')^{-1}$. Then it holds $|\hat{f}|=|g|<n$ (the inequality ``$<n$'' is valid due to Observation \ref{g kratsi}) and $|\hat{f}|_M=|f|_M-|f'|_M$.

Our goal is to estimate $|\hat{f}|_M-|g|_M$ which is equal to $|f|_M-|g|_M-|f'|_M$. For this purpose an estimate for $|f|_M-|g|_M$ will be needed; we obtain it from \eqref{*}: $|f|_M-|g|_M\geq(p+1)\cdot d+2-p$.

As for $|f'|_M$, since $f$ is followed by $S$, the last letter of $f'$ is $L$ (cf. Observation \ref{pripustne}).
With regard to this fact, Observation \ref{obsahM} implies $|f'|_M\leq\left\lceil\frac{d-1}{p^2+p}\right\rceil$. Now we distinguish two cases:\\
$\bullet$\quad If $d=1$, then $|f'|_M=0$, hence
$$
|\hat{f}|_M-|g|_M=|f|_M-|g|_M\geq(p+1)\cdot 1+2-p\geq3\,,
$$
$\bullet$\quad if $d\geq2$, then $|f'|_M\leq 1+\frac{d-1}{p^2+p}$, hence
\begin{multline*}
|\hat{f}|_M-|g|_M \geq |f|_M-|g|_M-\left(1+\frac{d-1}{p^2+p}\right) \geq (p+1)d+2-p-\left(1+\frac{d-1}{p^2+p}\right) =\\
=\left(p+1-\frac{1}{p^2+p}\right)(d-1)+2 \geq p+1-\frac{1}{p^2+p}+2 > p+2 \geq 4\,.
\end{multline*}
We see that for any value of $d>0$, the factors  $\hat{f}$ and $g$ of $u^{(p)}$ are of the same length less than $n$ and satisfy \eqref{unbalance M}. This is a contradiction with the minimality of $n$.\\

\pfk

\subsubsection*{Stage 2: Introduction of $x$, $y$}

At this moment we define another pair of factors.  Since $SfS$ and $SgS$ are factors of $u^{(p)}$, Observation \ref{JednozVzor} says that there exist factors $x$ and $y$ of $u^{(p)}$ such that
\begin{equation}\label{x}
\varphi_p(x)=fS
\end{equation}
and
\begin{equation}\label{y}
\varphi_p(y)=gS\,.
\end{equation}

\begin{prop}\label{x,y kratsi}
It holds $|x|\leq n-2(p^2-3)$ and $|y|\leq n-(2p^2-3)$, i.e. both factors $x$, $y$ are shorter that $v$, $w$.
\end{prop}
\pf
We use Propositions \ref{g kratsi} and \ref{f-g<=0}:
\begin{align*}
|y|&\leq|\varphi_p(y)|=|gS|=|g|+1\leq n-2(p^2-1)+1=n-(2p^2-3)\,, \\
|x|&\leq|\varphi_p(x)|=|fS|=|f|+1\leq|g|+1\leq n-(2p^2-3)\,.
\end{align*}
We remark that one can achieve better estimates, but these will be sufficient.
\pfk

\begin{prop}\label{vztah xy}
It holds
\begin{gather}
|x|_M-|y|_M=-(|f|-|g|)+|f|_M-|g|_M+p+1 \label{xyMfg}\\
|x|-|y|=1+|f|_M-|g|_M \label{xyfg}
\end{gather}
\end{prop}

\pf
The proof of both statements is straightforward using the definitions of $x$ and $y$, Proposition \ref{vzor}, the identity $|v|=|v|_L+|v|_S+|v|_M$ holding for any factor $v$, and Equality \eqref{f,g S}.
\pfk
Proposition \ref{vztah xy} has an immediate corollary:
\begin{prop}
\begin{gather}
|x|_M-|y|_M=|x|-|y|-(|f|-|g|)+p \label{xyMrov}\\
|x|_M-|y|_M\geq|x|-|y|+p \label{xyMner}
\end{gather}
\end{prop}
\pf
If we substitute for $|f|_M-|g|_M$ from \eqref{xyfg} to \eqref{xyMfg}, we obtain \eqref{xyMrov}. Equality \eqref{xyMrov} and Proposition \eqref{f-g<=0} then give \eqref{xyMner}.
\pfk

In what follows we split the proof according to the signum of $|x|-|y|$, and we show that whatever the value of $|x|-|y|$ is, it always contradicts the minimality of $n$.

\subsubsection*{The case $|x|-|y|>0$}

We set $|x|-|y|=d\geq1$ and denote $\hat{y}=yy'$, where $y'$ is a factor of $u^{(p)}$ of the length $d$ such that $yy'$ is a factor of $u^{(p)}$. Then $|x|=|\hat{y}|$. Moreover $|x|=|\hat{y}|<n$ due to Proposition \ref{x,y kratsi}.

Our goal is to estimate $|x|_M-|\hat{y}|_M$ which is equal to $|x|_M-|y|_M-|y'|_M$. The difference $|x|_M-|y|_M$ can be estimated using Inequality \eqref{xyMner}: $|x|_M-|y|_M\geq d+p$.

Let us proceed to $|y'|_M$. Recall at first that $\varphi_p(y)=gS$ (cf. \eqref{y}). Thus the last letter of the factor $y$ is either $L$ or $M$, which implies, with respect to Observation \ref{pripustne}, that the first letter of $y'$ is either $L$ or $S$, i.e. different from $M$. Observation \ref{obsahM} then gives $|y'|_M\leq \left\lceil\frac{d-1}{p^2+p}\right\rceil$. Now we distinguish two cases:\\
$\bullet$\quad If $d=1$, then $|y'|_M=0$, hence $|x|_M-|\hat{y}|_M=|x|_M-|y|_M\geq 1+p\geq3$;\\
$\bullet$\quad if $d\geq2$, then $|y'|_M\leq 1+\frac{d-1-1}{p^2+p}$, hence
\begin{multline*}
|x|_M-|\hat{y}|_M\geq d+p-1-\frac{d-2}{p^2+p}\geq p+(d-1)\left(1-\frac{1}{p^2+p}\right)+\frac{1}{p^2+p}\geq\\
\geq p+1\cdot\left(1-\frac{1}{p^2+p}\right)+\frac{1}{p^2+p}=p+1\geq3\,.
\end{multline*}
The factors $x$ and $\hat{y}$ are of the same length less than $n$ and for any $d>0$ satisfy $|x|_M-|\hat{y}|_M\geq 3$. In other words, they contradict the minimality of $n$.

\subsubsection*{The case $|x|-|y|=0$}

Equality \eqref{xyMrov} implies $|x|_M-|y|_M=-(|f|-|g|)+p$, and we know from Proposition \ref{f-g<=0} that $|f|-|g|\leq0$. Therefore:\\
$\bullet$\quad If $|f|-|g|\leq-1$ or $p\geq3$, we have $|x|_M-|y|_M\geq3$. Since moreover $|x|=|y|<n$, we have arrived at a contradiction with the minimality of $n$.\\
$\bullet$\quad If $|f|-|g|=0$ and $p=2$, Equality \eqref{xyfg} gives $|f|_M-|g|_M=-1$. This, however, contradicts Inequality \eqref{*}.

\bigskip

The last case to deal with is $|x|-|y|<0$. Since the situation $|x|-|y|=-1$ is very complicated to study, we will start with the case $|x|-|y|\leq-2$ and then deal with $|x|-|y|=-1$ separately.

\subsubsection*{The case $|x|-|y|\leq-2$}

\noindent Let us set for simplicity $d=|x|-|y|$.

Equality \eqref{xyfg} gives $|f|_M-|g|_M=d-1\leq-3$. We infer from here and from the minimality of $n$ that $|f|<|g|$ (otherwise one could consider the prefix of $f$ of the length $|g|$, denote it by $\check{f}$, then $|\check{f}|=|g|<n$, $|\check{f}|_M-|g|_M\leq|f|_M-|g|_M\leq-3$, which is a contradiction with the minimality of $n$).

Similarly, we infer from $|x|<y|<n$ and from the minimality of $n$ that \linebreak $|x|_M-|y|_M\leq2$ (otherwise we denote the prefix of $y$ of the length $|x|$ by $\check{y}$, then $|x|=|\check{y}|<n$, $|x|_M-|\check{y}|_M\geq 3$). Equality \eqref{xyMfg} then gives an upper bound on $|g|-|f|$, namely $|g|-|f|\leq2-p-d$.

Let us denote $\hat{f}=Sff'$, where $Sff'$ is a factor of $u^{(p)}$ and $|f'|=|g|-|f|-1$. Then $|f'|\geq0$ and $|\hat{f}|=|g|<n$, and it follows from \eqref{f} that if $f'\neq\epsilon$, the factor $f'$ can be chosen such that its first letter is $S$, i.e. different from $M$.

Now we are going to express $|\hat{f}|_M-|g|_M$ which is equal to $|f|_M-|g|_M+|f'|_M$. Let us distinguish the cases $|f'|\in\{0,1\}$ and $|f'|\geq2$.

$\bullet$\quad If $|f'|\in\{0,1\}$, we have $f'=\epsilon$ or $f'=S$, hence $|\hat{f}|_M-|g|_M=|f|_M-|g|_M=d-1\leq-3$\,.

$\bullet$\quad Let $|f'|\geq2$. Note that at the same time it holds $|f'|=|g|-|f|-1\leq2-p-d-1\leq -d-1$, thus necessarily $d\leq-3$. Observation \ref{obsahM} leads to
$$
|f'|_M\leq\left\lceil\frac{|f'|-1}{p^2+p}\right\rceil\leq\left\lceil\frac{2-p-d-1-1}{p^2+p}\right\rceil=\left\lceil\frac{-p-d}{p^2+p}\right\rceil\leq1+\frac{-p-d-1}{p^2+p}\,.
$$
Then
\begin{multline*}
|\hat{f}|_M-|g|_M= |f|_M-|g|_M+|f'|_M=d-1+|f'|_M\leq d-1+1+\frac{-p-d-1}{p^2+p}=\\
=d\cdot\frac{p^2+p-1}{p^2+p}-\frac{p+1}{p^2+p}\leq-3\cdot\frac{p^2+p-1}{p^2+p}-\frac{p+1}{p^2+p}=-3+\frac{2-p}{p^2+p}\leq-3\,.
\end{multline*}
We conclude that for any $d\leq-2$ it holds $|\hat{f}|=|g|<n$ and $|\hat{f}|_M-|g|_M\leq-3$. This is a contradiction with the minimality of $n$.

\subsubsection*{The case $|x|-|y|=-1$}

In this case, Equality \eqref{xyfg} implies $|f|_M-|g|_M=-2$ and Equality \eqref{xyMrov} implies $|x|_M-|y|_M=-1-(|f|-|g|)+p$.

Since $|x|<|y|<n$, necessarily $|x|_M-|y|_M\leq2$ (the contrary conradicts the minimality of $n$), hence \eqref{xyMrov} gives
\begin{equation}\label{-1fg}
|f|-|g|\geq p-3\,.
\end{equation}
We observe that this case, namely $|x|-|y|=-1$, can occur only for $p=2$, because:
\begin{itemize}
\item If $p$ was greater than 3, Inequality \eqref{-1fg} would contradict the inequality $|f|-|g|\leq0$, derived in Proposition \ref{f-g<=0}.
\item If $p$ was equal to 3, Inequality \eqref{-1fg} together with Proposition \ref{f-g<=0} would give $|f|-|g|=0$, which would not conform to \eqref{*}.
\end{itemize}
This allows us to restrict our considerations on the case $p=2$. Inequality \eqref{-1fg} implies $|f|-|g|\geq-1$, Inequality \eqref{*} implies $|f|-|g|\leq\frac{|f|_M-|g|_M}{3}=-\frac{2}{3}$, putting it together we infer $|f|-|g|=-1$. This allows us to use Eq. \eqref{xyMfg} to compute $|x|_M-|y|_M$: we obtain $|x|_M-|y|_M=2$.

Let us sum up the relations between $x$ and $y$ and between $f$ and $g$:
\begin{align}
|x|_M-|y|_M=2 \quad&\wedge\quad |x|-|y|=-1\,, \label{rozdxy}\\
|f|_M-|g|_M=-2 \quad&\wedge\quad |f|-|g|=-1\,. \label{rozdfg}
\end{align}

Our next goal is to prove this proposition:
\begin{prop}\label{slozite}
It holds:
\begin{itemize}
\item The word $Lx$ is a factor of $u^{(p)}$ and at it holds $Lx=LLS\cdots ML$,
\item the word $MLy$ is a factor of $u^{(p)}$ and at it holds $MLy=MLS\cdots SL$.
\end{itemize}
\end{prop}

\subsubsection*{Proof of Proposition \ref{slozite}.}
\noindent The proof will be done in ten steps.

\bigskip

\noindent \emph{Step 1.} (Possible prefixes and suffixes of $v$, $w'$)\\
The following four statements hold:
\begin{itemize}
\item[\textit{(i)}] \quad $v=MLSLLSM\cdots$ \qquad or \qquad $v=MLLSLLSM\cdots$,
\item[\textit{(ii)}] \quad $v=\cdots MLSLLSM$ \qquad or \qquad $v=\cdots MLLSLLSM$,
\item[\textit{(iii)}] \quad $w'=MLLSLLSM\cdots$ \qquad or \qquad $w'=MLSLLSLLSM\cdots$,
\item[\textit{(iv)}] \quad $w'=\cdots MLLSLLSM$ \qquad or \qquad $w'=\cdots MLSLLSLLSM$.
\end{itemize}

\pf
We have found in the proof of Observation \ref{g kratsi} that
$$
w'=Mw^{(1)}ww^{(2)}M=Mw^{(1)}\hat{w}^{(1)}M\cdots M\hat{w}^{(2)}w^{(2)}M\,,
$$
where $|w^{(j)}\hat{w}^{(j)}|_M=0$ and $|\hat{w}^{(j)}|\geq p^2+p$ for $j=1,2$. Taking into account Observation \ref{MzM}, we deduce that for $p=2$, the prefix of $w'$ and its suffix can be equal to either $MLLSLLSM$ or $MLSLLSLLSM$, thus (iii) and (iv) are proved.

Statements (i) and (ii) say in fact, with regard to Observation \ref{MzM}, that $v$ cannot have the segment $MLSLLSLLSM$ as its prefix and suffix, respectively. Let us suppose to the contrary that e.g. $v=MLSLLSLLSM\cdots$. We define $\hat{v}=(MLSLLSLLS)^{-1}v$, then $|\hat{v}|_M=|v|_M-1$. At the same time we denote $\hat{w}=\breve{w}^{-1}w$, where $\breve{w}$ is the prefix of $w$ of the length 9 (thus $|\hat{w}|=|\hat{v}|<n$). Since $w=\hat{w}^{(1)}M\cdots$, it follows from Observation \ref{MzM} that $\breve{w}$ contains exactly $1$ letter $M$, hence $|\hat{w}|_M=|w|_M-1$. Therefore $|\hat{v}|_M-|\hat{w}|_M=3$, which is a contradiction with the minimality of $n$.
\pfk

\noindent \emph{Step 2.} (Possible prefixes and suffixes of $f$, $g$)\\
It holds
\begin{itemize}
\item[\textit{(i)}] \quad $f=MLS\cdots$ \qquad or \qquad $f=LLS\cdots$,
\item[\textit{(ii)}] \quad $f=\cdots SML$ \qquad or \qquad $f=\cdots SLL$,
\item[\textit{(iii)}] \quad $g=LLS\cdots$ \qquad or \qquad $g=MLLS\cdots$,
\item[\textit{(iv)}] \quad $g=\cdots SLL$ \qquad or \qquad $g=\cdots SMLL$.
\end{itemize}

\pf
It is a trivial consequence of Step 1 and of the definitions of $f$ and $g$.
\pfk

\noindent \emph{Step 3.} (The prefix and suffix of $f$)\\
It holds $f=LLS\cdots SLL$.

\pf
We show that $f=LLS\cdots$, the proof of $f=\cdots SLL$ would be similar.

Step 2 implies that $f=MLS\cdots L$ or $f=LLS\cdots L$. Let us suppose for a while that $f=MLS\cdots L$. We put $\hat{f}=M^{-1}fS$ (it is a factor of $u^{(p)}$, because $fS$ is a factor of $u^{(p)}$, cf. \eqref{x}). Then $|\hat{f}|=|f|$, $|\hat{f}|_M=|f|_M-1$.

Since $g=\cdots LL$ according to Step 2, we are allowed to set $\hat{g}=gL^{-1}$. Then it holds $|\hat{g}|=|g|-1$, $|\hat{g}|_M=|g|_M$.

The factors $\hat{f}$ and $\hat{g}$ satisfy $|\hat{f}|=|\hat{g}|<n$ and $|\hat{f}|_M-|\hat{g}|_M=|f|_M-|g|_M-1=-3$, see \eqref{rozdfg}. This is a contradiction with the minimality of $n$.
\pfk

\noindent \emph{Step 4.} (The prefixes and suffixes of $v$, $w'$ and $g$)\\
It holds
\begin{itemize}
\item[\textit{(i)}] $v=MLLSLLSM\cdots MLLSLLSM$,
\item[\textit{(ii)}] $w'=MLSLLSLLSM\cdots MLSLLSLLSM$,
\item[\textit{(iii)}] $g=MLLS\cdots SMLL$.
\end{itemize}

\pf
Statement (i) is a straightforward consequence of Step 3. Statement (ii) follows from (i), which can be proven by contradiction using similar ideas as in the proof in Step 1. Statement (iii) is a consequence of (ii).
\pfk

\noindent \emph{Step 5.} (The prefixes and suffixes of $x$ and $y$)\\
It holds
\begin{itemize}
\item[\textit{(i)}] $x=L\cdots LL$ \qquad or \qquad $x=L\cdots ML$,
\item[\textit{(ii)}] $y=SLLS\cdots LSL$.
\end{itemize}

\pf
\textit{(i)} Since $fS=\varphi_2(x)$ and $fS=LLS\cdots SLLS$ according to Step 3, obviously $x=L\cdots LL$ or $x=L\cdots ML$, cf. Observation \eqref{pripustneVzor} for $p=2$.\\
\textit{(ii)} It follows from Step 4 that $y=SL\cdots SL$. Observation \ref{pripustne} implies $y=\cdots LSL$, Observation \ref{SzS} gives $y=SLLS\cdots$.
\pfk

\noindent \emph{Step 6.} (Two letters after $y$)\\
The word $yLS$ is a factor of $u^{(p)}$.

\pf
It is a consequence of Step 5 and Observation \ref{SzS}.
\pfk

\noindent \emph{Step 7.} (The suffix of $x$)\\
It holds $x=\cdots ML$.

\pf
We consider the result of Step 5 and prove that $x=\cdots LL$ contradicts the minimality of $n$. Let $x=\cdots LL$. Observation \ref{MzM} then implies that $x$ is followed either by $SM$ or by $SLLSM$. In the first case we define $\hat{x}=L^{-1}xSM$. Taking \eqref{rozdxy} into account, we see that the pair $\hat{x}$, $y$ contradicts the minimality of $n$. In the case when $x$ is followed by the group $SLLSM$, we define $\hat{x}=\breve{x}^{-1}xSLLSM$ where $\breve{x}$ is the prefix of $x$ of the length 2 (i.e. $\breve{x}=LL$ or $\breve{x}=LS$). Also we define $\hat{y}=yLS$ (this is allowed due to Step 6). The pair $\hat{x}$, $\hat{y}$ now contradicts the minimality of $n$ with regard to \eqref{rozdxy}.
\pfk

\noindent \emph{Step 8.} (The letter before $x$)\\
The word $Lx$ is a factor of $u^{(p)}$.

\pf
Since $SfS$ is a factor of $u^{(p)}$ and $fS=\varphi_p(x)$, it follows from Observation \ref{pripustneVzor} that $Lx$ or $Mx$ is a factor of $u^{(p)}$. However, if $Mx$ is a factor of $u^{(p)}$, then the factors $Mx$ and $y$ contradict the minimality of $n$ with regard to \eqref{rozdxy}.
\pfk

\noindent \emph{Step 9.} (The prefix of $x$)\\
It holds $x=LS\cdots$.

\pf
Since $Lx$ is a factor of $u^{(p)}$ according to Step 8 and $x=L\cdots$ according to Step 5, the word $Lx=LL\cdots$ is a factor of $u^{(p)}$. Observation \ref{pripustneVzor} then gives $Lx=LLS\cdots$.
\pfk

\noindent \emph{Step 10.} (Two letters before $y$)\\
The word $MLy$ is a factor of $u^{(p)}$.

\pf
Step 5 and Observation \ref{SzS} imply that $LLy$ or $MLy$ is a factor of $u^{(p)}$. Let us suppose for a while that $LLy=LLSLLS\cdots$ is a factor of $u^{(p)}$. Observation \ref{MzM} implies immediately that $LLy=LLSLLSM\cdots$. We introduce the word $\hat{y}=(SLLSM)^{-1}yLS$ which is a factor of $u^{(p)}$ due to Step 6. Then we define $\hat{x}=(LS)^{-1}x$, this is a factor of $u^{(p)}$ due to Step 9. It follows from \eqref{rozdxy} that the factors $\hat{x}$, $\hat{y}$ satisfy $|\hat{x}|=|\hat{y}|<n$, $|\hat{x}|_M-|\hat{y}|_M=3$, i.e. they contradict the minimality of $n$.
\pfk

This finishes the proof of Proposition \ref{slozite}. The statement follows from the results of Steps 5, 7, 8, 9, 10.


\subsubsection*{Stage 3: Introduction of $r$, $s$}

Proposition \ref{slozite} together with Observation \ref{JednozVzor} allow one to define factors $r$ and $s$ of $u^{(p)}$ such that
\begin{equation*}
\varphi_2(r)L=Lx\,, \qquad \varphi_2(s)L=Ly\,.
\end{equation*}
It is obvious that $|r|\leq|x|$, $|s|\leq|y|$, hence $|r|<n$, $|s|<n$, cf. Proposition \ref{x,y kratsi}.

Proposition \ref{vzor} and relations \eqref{rozdfg}, \eqref{rozdxy} enable us to compute $|r|-|s|$,
\begin{multline*}
|r|-|s|=(|x|_S+|x|_M)-(|y|_S+|y|_M)=(|x|_S-|y|_S)+(|x|_M-|y|_M)=\\
=(|fS|_M-|gS|_M)+(|x|_M-|y|_M)=(|f|_M-|g|_M)+(|x|_M-|y|_M)=-2+2=0\,,
\end{multline*}
and also $|r|_M-|s|_M$:
\begin{multline*}
|r|_M-|s|_M=(-|x|_L+2|x|_S)-(-|y|_L+2|y|_S)=-(|x|_L-|y|_L)+2(|x|_S-|y|_S)=\\
=-\left[|x|-|y|-(|x|_S-|y|_S)-(|x|_M-|y|_M)\right]+2(|x|_S-|y|_S)=\\
=-(|x|-|y|)+3(|x|_S-|y|_S)+(|x|_M-|y|_M)=\\
=-(|x|-|y|)+3(|fS|_M-|gS|_M)+(|x|_M-|y|_M)=-(-1)+3\cdot(-2)+2=-3\,.
\end{multline*}
The pair $r$, $s$ thus satisfy $|r|=|s|<n$ and $|r|_M-|s|_M<-3$. This is a contradiction with the minimality of $n$.

The proof of Theorem \ref{balanceM} is completed.


\section{Balance bound with respect to the letter $S$}

Once we know that $u^{(p)}$ is $2$-balanced with respect to the letter $M$, it is easy to prove that it is $2$-balanced with respect to the letter $S$ as well.

\begin{thm}\label{balanceS}
Let $v$, $w$ be factors of $u^{(p)}$ such that $|v|=|w|$. Then
$$
\left|\, |v|_S-|w|_S\, \right|\leq2\,.
$$
\end{thm}

\pf

We will proceed by contradiction. Let us assume that there exist factors $|v|$, $|w|$ of $u^{(p)}$ such that $|v|=|w|$ and $|v|_S-|w|_S>2$. Obviously one can suppose that $|v|_S-|w|_S=3$. Then, with regard to the substitution rule \eqref{varphi}, one has
\begin{align*}
|\varphi_p(v)|&=(p+1)|v|_L+|v|_S+p|v|_M\\
|\varphi_p(w)|&=(p+1)|w|_L+|w|_S+p|w|_M\,,
\end{align*}
hence, using the identities $|v|_L=|v|-|v|_S-|v|_M$, $|w|_L=|w|-|w|_S-|w|_M$,
\begin{align*}
|\varphi_p(v)|-|\varphi_p(w)|=&(p+1)(|v|_L-|w|_L)+(|v|_S-|w|_S)+p(|v|_M-|w|_M)=\\
=&(p+1)(|v|-|w|)-p(|v|_S-|w|_S)-(|v|_M-|w|_M)\,.
\end{align*}
Since the factors $v$ and $w$ are of the same length, they satisfy $|v|_M-|w|_M\geq-2$ with regard to Theorem \ref{balanceM}, hence
$$
|\varphi_p(v)|-|\varphi_p(w)|\leq(p+1)\cdot0-p\cdot3-(-2)=2-3p<0\,.
$$
Furthermore, since $|\varphi_p(v)|_M=|v|_S$ and $|\varphi_p(w)|_M=|w|_S$, one has
$$
|\varphi_p(v)|_M-|\varphi_p(w)|_M=|v|_S-|w|_S=3\,.
$$
Let $\hat{w}'$ be a prefix of $\varphi_p(w)$ of the length $|\varphi_p(v)|$. Then $|\hat{w}'|_M\leq|\varphi_p(w)|_M$, hence
$$
|\varphi_p(v)|_M-|\hat{w}'|_M\geq|\varphi_p(v)|_M-|\varphi_p(w)|_M\geq3\,,
$$
and since it holds $|\hat{w}'|=|\varphi_p(v)|$, we have arrived at a contradiction with Theorem~\ref{balanceM}.

\pfk


\section{Balance bound with respect to the letter $L$}\label{balL}

This section is devoted to the proof of the third part of Theorem \ref{balance}. As we will see, we will use both facts proved in the previous two sections.

\begin{thm}\label{balanceL}
Let $v$, $w$ be factors of $u^{(p)}$ such that $|v|=|w|$. Then
$$
\left|\, |v|_L-|w|_L\, \right|\leq3\,.
$$
\end{thm}

\pf

We will again proceed by contradiction. Let us assume that there exist factors $|v|$, $|w|$ of $u^{(p)}$ such that $|v|=|w|=n$ and
\begin{equation}\label{unbalance L}
|v|_L-|w|_L>3\,.
\end{equation}
Let $n$ be minimal number with this property.

We denote $v=v_1\cdots v_n$, $w=w_1\cdots w_n$. The minimality of $n$ implies
\begin{gather}
v_1=v_n=L\,, \label{vLL} \\
w_1\neq L\,,\quad w_n\neq L\,, \label{wxLxL}\\
|v|_L-|w|_L=4\,. \label{vw3L}
\end{gather}

The identity $|v|-|w|=|v|_L-|w|_L+|v|_S-|w|_S+|v|_M-|w|_M$ gives
$$
|v|_S-|w|_S+|v|_M-|w|_M=-4\,,
$$
and since $|v|_S-|w|_S\geq-2$ and $|v|_M-|w|_M\geq-2$ by virtue of Theorems \ref{balanceM} and \ref{balanceS}, respectively, we infer
$$
|v|_S-|w|_S=-2 \qquad \text{and} \qquad |v|_M-|w|_M=-2\,.
$$

Note that since $v_1=v_n=L$, Observation \ref{pripustneVzor} implies that there are $k,\ell\in\{0,1,\ldots,p-1\}$ and $X\in\{S,M\}$ such that $XL^{k}vL^\ell S$ is a factor of $u^{(p)}$.
We employ these $k$ and $\ell$ and define factors $v'$ and $w'$ of $u^{(p)}$ in the following way:
\begin{align*}
v'&=L^kvL^\ell S\,, \\
w'&=w_2\cdots w_n\,.
\end{align*}
Now we can apply Observation \ref{JednozVzor} which establishes the existence of factors $x$ and $y$ of $u^{(p)}$ satisfying
$$
\varphi_p(x)=v'\,, \qquad \varphi_p(y)=w'\,;
$$
obviously $|y|\leq |w'|=n-1$.

We are going to compute $|x|_L-|y|_L$ and $|x|-|y|$. For that purpose the following relations will be useful:
\begin{equation*}
\begin{array}{ll}
|v'|_L=|v|_L+k+\ell\,, & \qquad |w'|_L=|w|_L\,, \\
|v'|_S=|v|_S+1\,, & \qquad |w'|_S-|w|_S=-|w_1|_S\in\{-1,0\}\,, \\
|v'|_M=|v|_M\,, & \qquad |w'|_S+|w'|_M=|w|_S+|w|_M-1\,.
\end{array}
\end{equation*}
Employing these relations and Propositon \ref{vzor}, one can derive
\begin{multline*}
|x|_L-|y|_L=|v'|_L-(p-1)|v'|_S-\left(|w'|_L-(p-1)|w'|_S\right)=\\
=|v'|_L-|w'|_L-(p-1)\left(|v'|_S-|w'|_S\right)=\\
=|v|_L+k+\ell-|w|_L-(p-1)\left(|v|_S+1-|w|_S+|w_1|_S\right)=\\
=4+k+\ell-(p-1)\left(-2+1+|w_1|_S\right)\geq
4+0+0-(p-1)(-2+1+1)=4
\end{multline*}
and
\begin{align*}
|x|-|y|&=|v'|_S+|v'|_M-\left(|w'|_S+|w'|_M\right)=\\
&=|v|_S+1+|v|_M-\left(|w|_S+|w|_M-1\right)=\\
&=|v|_S-|w|_S+|v|_M-|w|_M+2
=-2-2+2=-2\,.
\end{align*}

Since $x$ is shorter than $y$, we consider the prefix of the factor $y$ of the length $|x|$ and denote it by $\hat{y}$; it obviously holds $|\hat{y}|_L\leq|y|_L$. Therefore
$$
|x|=|\hat{y}|<n \quad \text{and} \quad |x|_L-|\hat{y}|_L\geq|x|_L-|y|_L=4\,;
$$
in other words, the factors $x$ and $\hat{y}$ contradict the minimality of $n$.

\pfk

\section{Optimality of the balance bounds}

To complete the proof of Theorem \ref{balance}, we have to show that the bounds 3, 2 and 2 corresponding to $L$, $S$ and $M$, respectively, are optimal. To demonstrate this fact it suffices to find three pairs of factors of $u^{(p)}$, let us denote them by $(v^{(M)},w^{(M)})$, $(v^{(S)},w^{(S)})$ and $(v^{(L)},w^{(L)})$, such that

\begin{equation*}
\begin{array}{ccc}
|v^{(M)}|=|w^{(M)}| \quad & \text{and} & \quad \left|\,|v^{(M)}|_M-|w^{(M)}|_M\,\right|=2\,, \\
|v^{(S)}|=|w^{(S)}| \quad & \text{and} & \quad \left|\,|v^{(S)}|_S-|w^{(S)}|_S\,\right|=2\,, \\
|v^{(L)}|=|w^{(L)}| \quad & \text{and} & \quad \left|\,|v^{(L)}|_L-|w^{(L)}|_L\,\right|=3\,.
\end{array}
\end{equation*}
There are many possibilities, one can take for example:
\begin{itemize}
\item $v^{(M)}=M\varphi_p^2(SM)$,\; $w^{(M)}=\varphi_p^2(SL)(L^{p-2}SM)^{-1}$\\
Then $|v^{(M)}|=|w^{(M)}|=p^2+p+1$,\; $|v^{(M)}|_M-|w^{(M)}|_M=2$.
\item $v^{(S)}=L^{-p}\varphi_p^2(LL)M^{-1}$,\; $w^{(S)}=\varphi_p^2(ML)L^p$\\
Then $|v^{(S)}|=|w^{(S)}|=2p^2+2p+1$,\; $|v^{(S)}|_S-|w^{(S)}|_S=2$.
\item $v^{(L)}=\varphi_p^3(ML^{p-1})\varphi_p^2(LL)(L^{p-2}SM)^{-1}$,\;
$w^{(L)}=SM\varphi_p^2(SM)\varphi_p^4(SM)$\\
Then $|v^{(L)}|=|w^{(L)}|=5p^2+6p+5$,\; $|v^{(L)}|_L-|w^{(L)}|_L=3$.
\end{itemize}


\section{Abelian complexity}\label{ACsect}

In this section we will determine the optimal bound for the Abelian complexity function $\mathrm{AC}(n)$ of the infinite word $u^{(p)}$, taking advantage of the results on optimal balance bounds that have been derived in the previous part of the paper.

The set of factors of $u^{(p)}$ of the length $n$ will be denoted by $\mathcal{F}_u(n)$, and the symbol $\mathcal{P}_u(n)$ will stand for the set of corresponding Parikh vectors, i.e. \linebreak $\mathcal{P}_u(n)=\left\{\Psi(w)\,\left|\,w\in\mathcal{F}_u(n)\right.\,\right\}$.

\begin{prop}\label{mozPar}
For each $n\in\N$ there are numbers $s_n\in\N$ and $m_n\in\N$ such that $\mathcal{P}_u(n)\subset\{\Psi^{(n)}_1,\ldots,\Psi^{(n)}_9\}$, where
\begin{equation}\label{mozneParikh}
\begin{array}{rlcccccl}
\Psi^{(n)}_1=&(& n-s_n-m_n+2 &,& s_n-1 &,& m_n-1 &)\,, \\
\Psi^{(n)}_2=&(& n-s_n-m_n+1 &,& s_n-1 &,& m_n &)\,, \\
\Psi^{(n)}_3=&(& n-s_n-m_n &,& s_n-1 &,& m_n+1 &)\,, \\
\Psi^{(n)}_4=&(& n-s_n-m_n+1 &,& s_n &,& m_n-1 &)\,, \\
\Psi^{(n)}_5=&(& n-s_n-m_n &,& s_n &,& m_n &)\,, \\
\Psi^{(n)}_6=&(& n-s_n-m_n-1 &,& s_n &,& m_n+1 &)\,, \\
\Psi^{(n)}_7=&(& n-s_n-m_n &,& s_n+1 &,& m_n-1 &)\,, \\
\Psi^{(n)}_8=&(& n-s_n-m_n-1 &,& s_n+1 &,& m_n &)\,, \\
\Psi^{(n)}_9=&(& n-s_n-m_n-2 &,& s_n+1 &,& m_n+1 &)\,.
\end{array}
\end{equation}
Consequently, $\mathrm{AC}(n)\leq9$.
\end{prop}

\pf
For each $n\in\N$, Theorem \ref{balanceS} implies that
$$
\max\{\left.|v|_S\,\right|\,v\in\mathcal{F}_u(n)\}-\min\{\left.|v|_S\,\right|\,v\in\mathcal{F}_u(n)\}\leq2\,,
$$
thus there is an $s_n\in\N$ such that
$$
v\in\mathcal{F}_u(n) \quad \Rightarrow \quad |v|_S\in\{s_n-1,s_n,s_n+1\}\,.
$$
Similarly, using Theorem \ref{balanceM}, one finds that there is an $m_n\in\N$ such that
$$
v\in\mathcal{F}_u(n) \quad \Rightarrow \quad |v|_M\in\{m_n-1,m_n,m_n+1\}\,.
$$
Since $|v|_L=|v|-|v|_S-|v|_M$, we deduce that the set $\mathcal{P}_u(n)$ is a subset of $\{\Psi^{(n)}_1,\ldots,\Psi^{(n)}_9\}$, where $\Psi^{(n)}_j$, $j=1,\ldots,9$, are given by \eqref{mozneParikh}.
\pfk

Proposition \ref{mozPar} gives an upper bound of $\mathrm{AC}(n)$, namely $\mathrm{AC}(n)\leq9$ for every $n\in\N$, but we can even say more. Indeed, Theorem \ref{balanceL} implies that $\mathcal{P}_u(n)$ cannot contain at the same time $\Psi^{(n)}_1$ and $\Psi^{(n)}_9$, hence $\mathrm{AC}(n)\leq8$. In fact, the optimal bound is even lower, as we will show with the help of the following proposition.

\begin{prop}\label{vylucujici}
There is no pair of factors $v$, $w$ of $u^{(p)}$ such that their Parikh vectors satisfy $\Psi(v)-\Psi(w)=(3,-2,-1)$.
\end{prop}

\pf
We prove the statement by contradiction. Let us suppose that there are factors $v$ and $w$ of $u^{(p)}$ satisfying
\begin{equation}\label{ACLSM}
|v|_L-|w|_L=3\,,\quad |v|_S-|w|_S=-2\,,\quad |v|_M-|w|_M=-1\,.
\end{equation}
The factors $v$ and $w$ are obviously of the same length which we denote by $n$. If $v=v_1\cdots v_n$ and $w=w_1\cdots w_n$, we may assume without loss of generality that $\{v_1,v_n\}\cap\{w_1,w_n\}=\emptyset$ (if e.g. $v_1=w_n$, we can replace the pair  $(v,w)$ by $(v_2\cdots v_n,w_1,\cdots w_{n-1})$). Having this assumption and Eqs. \eqref{ACLSM}, we deduce
\begin{equation}\label{neL}
w_1\neq L\quad \text{and} \quad w_n\neq L\,,
\end{equation}
because if e.g. $w_1=L$, then the pair $v_2\cdots v_n$ and $w_2\cdots w_n$ contradicts Theorem \ref{balanceL}.
Similarly, Theorem \ref{balanceS} implies that
\begin{equation}\label{neS}
v_1\neq S\quad \text{and} \quad v_n\neq S\,.
\end{equation}

Let $k,\ell\geq0$ and $X,Y\in\{S,M\}$ be such numbers and letters that the word $XL^kvL^\ell Y$ is a factor of $u^{(p)}$. It is a trivial fact that such $k,\ell$, $X,Y$ exist. Moreover, Observation \ref{pripustne} implies $Y=S$, Observation \ref{pripustneVzor} gives $k,\ell\leq p$. We take these $k,\ell$ and $X$ and define a factor $v'$ of $u^{(p)}$, and then define also a factor $w'$ of $u^{(p)}$:
\begin{gather*}
v'=L^k vL^\ell S\,, \\
w'=w_2\cdots w_n\,.
\end{gather*}
Now we apply Observation \ref{JednozVzor} which says that there are factors $x$ and $y$ of $u^{(p)}$ such that
$$
v'=\varphi_p(x)\,, \qquad w'=\varphi_p(y)\,.
$$
Let us write down the following relations between $v'$ and $v$ and between $w'$ and $w$,
\begin{equation*}
\begin{array}{ll}
|v'|_L=|v|_L+k+\ell\,, & \qquad |w'|_L=|w|_L\,, \\
|v'|_S=|v|_S+1\,, & \qquad |w'|_S-|w|_S=-|w_1|_S=|w_1|_M-1\in\{-1,0\}\,, \\
|v'|_M=|v|_M\,, & \qquad |w'|_S+|w'|_M=|w|_S+|w|_M-1\,,
\end{array}
\end{equation*}
and use them to express $|x|_M-|y|_M$,
\begin{equation}\label{ACpocetM}
\begin{split}
|x|_M-|y|_M&=-|v'|_L+p|v'|_S-\left(-|w'|_L+p|w'|_S\right)=\\
&=-\left(|v'|_L-|w'|_L\right)+p\left(|v'|_S-|w'|_S\right)=\\
&=-(|v|_L+k+\ell-|w|_L)+p\left(|v|_S+1-|w|_S+1-|w_1|_M\right)=\\
&=-(3+k+\ell)+p\left(-2+2-|w_1|_M\right)=-3-k-\ell-p|w_1|_M\,,
\end{split}
\end{equation}
(note that $|x|_M-|y|_M\leq-3$), and also to express $|x|-|y|$:
\begin{align*}
|x|-|y|&=|v'|_S+|v'|_M-\left(|w'|_S+|w'|_M\right)=\\
&=|v|_S+1+|v|_M-\left(|w|_S+|w|_M-1\right)=\\
&=|v|_S-|w|_S+|v|_M-|w|_M+2=-2-1+2=-1\,.
\end{align*}

\noindent\textit{Statement.} It holds
\begin{itemize}
\item[\textit{(i)}] $y=M\cdots M$,
\item[\textit{(ii)}] $w_1=S$.
\end{itemize}
Both (i) and (ii) can be proved by contradiction:

\textit{(i)} Suppose that e.g. the first letter of $y$ is $X\neq M$. We define $\hat{y}=X^{-1}y$, then $|x|=|\hat{y}|$ and $|x|_M-|\hat{y}|_M=|x|_M-|y|_M\leq-3$, which is a contradiction with Theorem \ref{balanceM}. Similarly we prove that the last letter of $y$ is $M$.

\textit{(ii)} Suppose that $w_1\neq S$. We know from \eqref{neL} that also $w_1\neq L$, hence $w_1=M$. Then \eqref{ACpocetM} gives $|x|_M-|y|_M\leq-3-p\leq-5$. Let $Y$ be the last letter of $y$ and $\hat{y}=yY^{-1}$. Then it holds $|\hat{y}|=|x|$ and $|x|_M-|\hat{y}|_M=|x|_M-|y|_M+|Y|_M\leq-5+1=-4$, which is a contradiction with Theorem \ref{balanceM}.

\bigskip

\noindent Statements (i) and (ii) can be used to determine the first three letters of $w$:
$$
w=w_1\varphi_p(y)=S\varphi_p(M\cdots M)=SL^{p-1}S\cdots L^{p-1}S\,.
$$
However, Observation \ref{SzS} implies that no factor of $u^{(p)}$ can contain the segment $SL^{p-1}S$. This is a contradiction, thus the proposition is proved.

\pfk

\begin{thm}
\textit{(i)} Let $u^{(p)}$ be the fixed point of the substitution $\varphi_p$ defined in \eqref{varphi}. Then its Abelian complexity satisfies
$\mathrm{AC}(n)\leq7$ for all $n\in\N$ and for all $p>1$.\\
\textit{(ii)} The bound $7$ is optimal, i.e. it cannot be improved. Moreover, for any $p>1$ there is infinitely many numbers $n\in\N$ such that $\mathrm{AC}(n)=7$.
\end{thm}

\pf
\textit{(i)} Since $\Psi^{(n)}_1-\Psi^{(n)}_8=(3,-2,-1)$ and $\Psi^{(n)}_2-\Psi^{(n)}_9=(3,-2,-1)$, Proposition \ref{vylucujici} implies that the set $\mathcal{P}_u(n)$ can contain at most one of the vectors $\Psi^{(n)}_1$, $\Psi^{(n)}_8$ and at most one of the vectors $\Psi^{(n)}_2$, $\Psi^{(n)}_9$. Therefore $\mathcal{P}_u(n)$ has at most 7 elements, hence $\mathrm{AC}(n)\leq7$.

\textit{(ii)} We prove (ii) by showing that there are infinitely many values $n\in\N$ with this property: There are $7$ factors of $u^{(p)}$ of the length $n$ such that their Parikh vectors are mutually different.

\noindent For all $N\in\N$, we define two auxiliary words:
\begin{gather*}
v^{(N)}=\varphi_p^{2N+2}(L)\varphi_p^{2N+1}(L)\,,\\
w^{(N)}=\left(\varphi_p^{2N+1}(L)\right)^{-1}\varphi_p^{2N+2}(LS)\varphi_p^{2N+1}(L)\varphi_p^{2N}(L)\,.
\end{gather*}
\textit{Statement.} For each $N\in\N$ it holds:

\textit{(i)} $v^{(N)}$ and $L^pSMv^{(N)}$ are factors of $u^{(p)}$,

\textit{(ii)} $w^{(N)}$ and $SML^{p-1}Sw^{(N)}$ are factors of $u^{(p)}$,

\textit{(iii)} $v^{(N)}$ has the suffix $SML^{p-1}S$ and $w^{(N)}$ has the suffix $L^pSM$,

\textit{(iv)} $\Psi(v^{(N)})=\Psi(w^{(N)})$.

To see (i), we observe that $v^{(N)}$ is a prefix of $\varphi_p^{2N+2}(LL)$, thus obviously a factor of $u^{(p)}$. Moreover, since $MLL$ is a factor of $u^{(p)}$ (cf. Observation \ref{MzM}), the word $\varphi_p^{2N+2}(MLL)$ is a factor of $u^{(p)}$ as well. Then it suffices to show that $\varphi_p^{2N+2}(M)$ has the suffix $L^pSM$ for each $N\in\N$, which can be done easily by induction.

The proof of (ii) is similar: Since $\varphi_p^{2N+1}(L)w^{(N)}$ is a factor of $\varphi_p^{2N+2}(LSL)$, $w^{(N)}$ is a factor of $u^{(p)}$.
To demonstrate that $SML^{p-1}Sw^{(N)}$ is a factor of $u^{(p)}$, it suffices to realize that $SML^{p-1}S$ is a suffix of $\varphi_p^{2N+1}(L)$ for all $N\in\N$, which can be shown by induction.

Let us proceed to (iii). By the definitions above, $v^{(N)}$ has the suffix $\varphi_p^{2N+1}(L)$, $w^{(N)}$ has the suffix $\varphi_p^{2N}(L)$. Now we can use the facts known from (i) and (ii), namely that $\varphi_p^{2N+1}(L)=\cdots SML^{p-1}S$ and $\varphi_p^{2N}(L)=\cdots L^pSM$ for all $N\in\N$.

The statement (iv) is a consequence of these two equalities:
\begin{equation*}
\begin{split}
v^{(N)}&=\varphi_p^{2N+2}(L)\varphi_p^{N+1}(L)=\varphi_p^{2N+2}(L)\varphi_p^{2N}(L^pS)=\varphi_p^{2N+2}(L)\left(\varphi_p^{2N}(L)\right)^p\varphi_p^{2N}(S)\,,\\
w^{(N)}&=\left(\varphi_p^{2N+1}(L)\right)^{-1}\varphi_p^{2N+2}(L)\varphi_p^{2N}(L^{p-1}S)\varphi_p^{2N+1}(L)\varphi_p^{2N}(L)=\\
&=\left(\varphi_p^{2N+1}(L)\right)^{-1}\varphi_p^{2N+2}(L)\left(\varphi_p^{2N}(L)\right)^{p-1}\varphi_p^{2N}(S)\varphi_p^{2N+1}(L)\varphi_p^{2N}(L)\,.
\end{split}
\end{equation*}

For each $N\in\N$, the auxiliary factors $v^{(N)}$ and $w^{(N)}$ are intrumental to define another set of words that we denote $f^{(1)},\ldots,f^{(7)}$:
\begin{gather*}
f^{(1)}=L^pSMv^{(N)}(SML^{p-1}S)^{-1}\,, \qquad f^{(2)}=Mv^{(N)}S^{-1}\,, \\
f^{(3)}=LSw^{(N)}(SM)^{-1}\,, \qquad f^{(4)}=v^{(N)}\,, \qquad f^{(5)}=SMv^{(N)}(LS)^{-1}\,, \\
f^{(6)}=Sw^{(N)}M^{-1}\,, \qquad f^{(7)}=SML^{p-1}Sw^{(N)}\left(L^pSM\right)^{-1}\,.
\end{gather*}

It follows from Statement above that all the words $f^{(1)},\ldots,f^{(7)}$ are factors of $u^{(p)}$ of the same length $n_N=|v^{(N)}|=|w^{(N)}|$. If we compare the numbers of letters $L$, $S$ and $M$ in these factors, we find that their Parikh vectors and mutually different. This means that $\mathrm{AC}(n)=7$ for every $n=|v^{(N)}|$, where $N\in\N$, i.e. the function $\mathrm{AC}$ attains the value $7$ infinitely many times.

\pfk

We finish this section by two statements describing the lower bound and the range of $\mathrm{AC}(n)$.

\begin{prop}
It holds $\mathrm{AC}(n)\geq 3$ for all $n\in\N$ and for all $p>1$.
\end{prop}

\pf
At first we show that for any prefix $v$ of $u^{(p)}$, both $LSv$ and $SMv$ are factors of $u^{(p)}$. This will be obvious from these three simple facts:\\
\textit{(i)} Since $LL$ is a factor of $u^{(p)}$, both $\varphi_p^{2N}(L)\varphi_p^{2N}(L)$ and $\varphi_p^{2N+1}(L)\varphi_p^{2N+1}(L)$ are factors of $u^{(p)}$ for any $N\in\N$.\\
\textit{(ii)} It holds $\varphi_p^{2N}(L)=\cdots SM$ and $\varphi_p^{2N+1}(L)=\cdots LS$ for any $N\in\N$.\\
\textit{(iii)} There is an $N\in\N$ such that $v$ is a prefix of $\varphi_p^{2N}(L)$, and thus of $\varphi_p^{2N+1}(L)$.

Statements (i) and (ii) imply that both $SM\varphi_p^{2N}(L)v$ and $LS\varphi_p^{2N+1}(L)v$ are factors of $u^{(p)}$, and (iii) then implies that $SMv$ and $LSv$ are factors of $u^{(p)}$ as well.

Let now $v=u_0u_1\cdots u_{n-1}$ be the prefix of $u^{(p)}$ of the length $n$ and let $\Psi(v)$ be its Parikh vector. We will show that then there are factors $v'$ and $v''$ of $u^{(p)}$ such that $\Psi(v')\neq\Psi(v)\neq\Psi(v'')$ and $\Psi(v')\neq\Psi(v'')$. We distinguish three cases according to $u_{n-1}$:\\
$\bullet$ If $u_{n-1}=L$, then we set $v'=Su_0u_1\cdots u_{n-2}$, $v''=Mu_0u_1\cdots u_{n-2}$.\\
$\bullet$ If $u_{n-1}=S$, it holds necessarily $u_{n-2}=L$ (cf. Observation \ref{pripustne}). In this case we set $v'=Mu_0u_1\cdots u_{n-2}$, $v''=SMu_0u_1\cdots u_{n-3}$.\\
$\bullet$ If $u_{n-1}=M$, it holds necessarily $u_{n-2}=S$ (cf. again Observation \ref{pripustne}). We set $v'=Su_0u_1\cdots u_{n-2}$, $v''=LSu_0u_1\cdots u_{n-3}$.
\pfk

\begin{pozn}
It can be demonstrated that for each value $k\in\{3,4,5,6,7\}$ there is an $n\in\N$ such that $\mathrm{AC}(n)=k$, but we omit the proof with regard to the length of the paper.
\end{pozn}

\section{Conclusion}

We have studied balance properties and the Abelian complexity of a certain class of infinite ternary words. We have found the optimal $c$ such that these words are $c$-balanced, and also the optimal bound for their Abelian complexity functions. We have introduced a new notion, namely the property that a word is ``$c$-balanced with respect to the letter $a$'', which helped us to proceed more effectively from the knowledge of the balance properties to the estimate on $\mathrm{AC}(n)$.

The class of words studied in this paper has one parameter $p>1$. However, it emerged from our calculations that all the three optimal bounds for balances with respect to particular letters, as well as the optimal bound for the Abelian complexity, are independent of the value of $p$.

The problem has one more aspect. So far the subword complexity and the Abelian complexity are considered as highly independent of each other (cf. e.g. the work of Richome, Saari and Zamboni). However, our result can indicate that there are connections between them, for the present waiting for their discovery. It has been recently shown in \cite{KP} that a fixed point of a canonical substitution associated with a non-simple cubic Parry number has affine factor complexity if and only if it belongs just to the class with which we have dealt in this paper. Therefore, briefly speaking, ``if the factor complexity is affine, then the Abelian complexity has the optimal bound $7$'' holds in the cubic non-simple Parry case. We remark that this sort of statement holds as well in the quadratic non-simple Parry case, although there it is a trivial fact.

\section*{Acknowledgements}

The author is grateful to K. B\v rinda for performing a numerical experiment and to E. Pelantov\'a for helpful discussions on various aspects of the problem.


\end{document}